# Remarques sur l'expression de la généralité en mathématiques

Alain Herreman[1]

Université Rennes 1 – CNRS

alain.herreman@univ-rennes1.fr

*Résumé*

Cet article s'attache à dégager sur quelques exemples les conditions de possibilité de l'expression de la généralité en mathématiques. Au cours du XXe siècle les mathématiciens on pu énoncer des théorèmes portant sur la totalité des théorèmes (théorèmes d'incomplétude de Gödel, etc.). L'introduction du théorème de Löwenheim-Skolem nous donnera l'occasion de saisir la mise en place des conditions de tels énoncés généraux. Avec les *Disquisitiones arithmeticae* nous verrons au contraire comment l'impossibilité d'utiliser le mode de représentation qui lui était familier conduisit Gauss à recourir à des ensembles.

## Introduction

*On dit qu'à force d'ascèse certains bouddhistes parviennent à voir tout un paysage dans une fève. C'est ce qu'auraient bien voulu les premiers analystes du récit : voir tous les récits du monde (il y en a tant et tant eu) dans une seule structure.* On comprend à partir de ces deux premières phrases de *S/Z* que Roland Barthes n'était pas de ceux qui croient qu'il est possible de voir tout un paysage dans une fève. Il ne nous dit pas si la fève fait partie du paysage ou si une même fève permet de voir tous les paysages. La fève du moine Bouddhiste illustre le pouvoir des représentations qui rendent présent, visible et rassemblent ce qui est autrement absent, invisible et dispersé. Mais son usage effectif est douteux et requiert sans doute une longue ascèse… Elle n'est guère susceptible d'être partagée ou transmise. C'est néanmoins le pouvoir de telles représentations que nous voulons considérer ici à partir du cas des mathématiques.

Nombre d'énoncés mathématiques sont à bien des égards remarquables en raison notamment de leur *généralité*. Un énoncé aussi simple que « toute figure rectiligne peut être transformée en un carré de même aire » (paraphrase de *Euclide* II-14) est un énoncé typique, à la fois commun en mathématiques et sans équivalent ailleurs. Sa généralité est bien différente par exemple de celle d'un énoncé comme « tous les hommes sont mortels ». Celui-ci renvoie à une induction fondée sur le constat selon lequel jusqu'à présent tous les hommes ont bien dû être mortels, celui-là sur une démonstration. Mais si leur différence ressort des justifications qui peuvent être données de chacun, celles-ci ne créent pas tant cette différence qu'elles ne déploient autrement les caractéristiques de chaque énoncé.

Nous nous proposons dans ce qui suit de dégager une condition qui rend compte de la généralité d'une large variété d'énoncés mathématiques.

## *I. L'expression de la généralité et le théorème de Löwenheim-Skolem*

Le théorème de Löwenheim-Skolem est un exemple de ces énoncés que l'on reconnaît comme étant mathématiques. En voici un énoncé :

> Si une théorie du premier ordre ayant un nombre au plus dénombrable d'axiomes admet un modèle, alors elle admet un modèle dénombrable.

C'est en tant que tel qu'il nous intéresse ici : un énoncé à la fois d'une généralité étonnante et tout à fait habituelle en mathématiques. Celui-ci porte sur *toutes* les théories (du premier ordre), d'autres portent sur *tous* les nombres, sur *toutes* les fonctions, etc. Les mathématiques ont donc une notion, ici, de *théorie* qui est telle qu'il soit possible de formuler à leur propos des énoncés généraux, mais aussi de les démontrer et en plus de les appliquer à la théorie des nombres entiers, à la théorie des nombres réels, à la théorie des groupes, à la théorie des ensembles etc. C'est une condition assez générale de la possibilité de ces énoncés que nous voulons dégager. Celle-ci une fois repérée, il sera ensuite facile de voir qu'elle est indépendante de l'énoncé choisi, et en particulier de son appartenance à la logique. Cette condition peut être dégagée aussi bien à partir des diverses démonstrations qui ont été données du théorème que de ses applications. Pour éviter les difficultés auxquelles conduirait l'analyse d'une de ces démonstrations, nous dégagerons cette condition à partir d'une des principales applications de ce théorème, ce qui sera l'occasion de rappeler le « paradoxe de Skolem » et de revenir sur son interprétation.

### Le théorème de Löwenheim-Skolem

Le théorème de Löwenheim-Skolem s'applique à une théorie comme celle de la théorie des nombres réels ou même de la théorie des ensembles. Il suppose que celle-ci admet un modèle ce qui est *a priori* le cas puisque la théorie aura généralement été développée parce que l'on en a un modèle, au moins intuitif. Il n'est en général pas dénombrable, c'est-à-dire qu'il contient plus d'éléments qu'il n'y a de nombres entiers[2]. Ainsi, puisqu'il y a plus de nombres réels que d'entiers, notre modèle des nombres réels ne sera pas dénombrable. *A fortiori*, un modèle de la théorie des ensembles, qui doit contenir tous les nombres réels, ne sera pas non plus dénombrable. Et pourtant, la conclusion du théorème est qu'*il existe alors un modèle dénombrable de cette théorie,* c'est-à-dire qu'il est possible d'interpréter la théorie de telle sorte qu'il n'y ait pas à considérer plus qu'un nombre dénombrable d'entités. Skolem précisera ce résultat et démontrera que les nombres entiers eux-mêmes peuvent servir d'interprétation, c'est-à-dire que l'on peut remplacer tous les objets mathématiques (nombres réels, fonctions, ensembles, etc.) par des nombres entiers. Un tel énoncé est concevable, énonçable, démontrable et, on le voit, applicable.

---
2   i.e. qu'il n'y a pas de bijection entre les éléments du modèle et l'ensemble des entiers naturels.

Ce théorème a été énoncé et démontré pour la première fois par Löwenheim dans un article publié en 1915 (« Über Möglichkeiten im Relativkalkül », *Mathematische Annalen* 76, 447-470). Löwenheim utilisait le formalisme logique développé par Schröder dans son *Algebra der Logik* (1890-1905). Sa démonstration fait intervenir des développements en sommes de produits de formules, au moyen d'indices d'indices, assortis d'un jeu de recodage assez obscur. En 1920, Skolem généralise le théorème et en donne une nouvelle démonstration, toujours dans la logique algébrique, mais en remplaçant les développements des produits de Löwenheim par des opérations ensemblistes introduites par Dedekind (*Was sind und was sollen* die *Zahlen?*, Braunschweig, 1888) et en recourant à l'axiome du choix. Deux ans plus tard il publie un autre article avec une démonstration simplifiée et ne recourant pas à l'axiome du choix. C'est là qu'il propose l'application que nous voulons considérer et qu'il expose son paradoxe.

Cet article, publié en 1922, est un texte éminemment rhétorique. Son propos est de convaincre les mathématiciens que les mathématiques ne sauraient être fondées sur la théorie des ensembles comme ils semblent généralement enclins à le croire. En effet, le développement de la théorie des ensembles au cours de la seconde moitié du 19ème siècle s'est aussi accompagné de la découverte d'antinomies. Henri Poincaré pouvait ainsi s'exclamer en 1906 : « *la logique n'est plus stérile, elle engendre l'antinomie* »[3]. En 1908, Ernst Zermelo, un élève de David Hilbert, proposa un système d'axiomes *ad hoc* visant à garantir les principales opérations ensemblistes assorties de quelques restrictions pour empêcher la reproduction des antinomies connues. Cette démarche fut semble-t-il couronnée de succès puisqu'en 1922 Skolem constate que les mathématiciens sont devenus nombreux à souscrire à l'idée que la théorie des ensembles offre un fondement satisfaisant des mathématiques. Son article est intégralement consacré à réfuter cette idée au moyen de cinq arguments. Le « paradoxe de Skolem » est l'un d'eux, et c'est pour pouvoir appliquer le théorème de Löwenheim dans ce cadre rhétorique particulier qu'il en donne une deuxième démonstration. L'intérêt de cet argument sur les quatre autres est de viser directement le système d'axiomes de Zermelo. Présentons donc ce paradoxe et l'argument associé. L'application du théorème, sous la forme généralisée par Skolem, aux axiomes de Zermelo établit l'existence d'un modèle dénombrable de toutes les mathématiques. Dans ce modèle, l'ensemble des nombres réels serait dénombrable. Un des premiers théorèmes établis par Cantor, à partir duquel peut se déployer la hiérarchie des cardinaux transfinis, démontre pourtant le contraire. Il s'agit bien seulement d'un paradoxe et non d'une contradiction. En effet, la propriété « être dénombrable » signifie ici qu'il existe une application bijective entre l'ensemble considéré, par exemple celui des nombres réels et l'ensemble des nombres entiers. Pour qu'il y ait une contradiction il faudrait encore que cette application, conçue comme un ensemble, fasse elle-même partie du modèle. Ce n'est en tout cas pas le cas de celle construite par la démonstration du théorème. Autrement dit, de l' « intérieur » du modèle on « ne voit pas » que l'ensemble des nombres réels du modèle est dénombrable, pas plus d'ailleurs qu'on ne voit que le modèle est lui-même dénombrable (il ne peut pas exister d'application *dans* le modèle ayant comme ensemble de départ ou d'arrivée le modèle lui-même). Ni de l' « intérieur », ni non plus de l' « extérieur », il n'y a de contradiction. Skolem ne fait pas passer ce paradoxe, qu'il présente comme tel et qu'il explique parfaitement, pour une contradiction. Son argument est rhétorique, mais pas

---
[3]  Poincaré, Henri, « Les mathématiques et la logique », *Revue de Métaphysique et de Morale* 14, 1906, 316 reproduit *in* Heinzmann, Gerhard, *Poincaré, Russell, Zermelo et Peano. Textes de la discussion (1906-1912) sur les fondements des mathématiques : des antinomies à la prédicativité*. Paris, Blanchard, 1986, 103,

sophistique! L'argument qu'il en tire contre l'idée que les axiomes de Zermelo fonderaient les mathématiques consiste à observer que *toutes les notions mathématiques deviendraient de ce fait relatives* ; car si l'on accepte cette hypothèse, il existe un *deuxième* modèle des mathématiques et même bien d'autres (seule la pluralité importe pour l'argument, et non plus la cardinalité du modèle). Plus aucun objet mathématique n'a de référence *unique,* toutes les mathématiques deviennent plurivoques. Il serait dès lors établi, avec la certitude d'un théorème, que les mathématiciens ne savent pas de quoi ils parlent, qu'ils ne parlent pas nécessairement de la même chose, et que le savoir n'y change rien. En appliquant son théorème aux axiomes de Zermelo, Skolem peut établir qu'accepter ces axiomes comme fondement des mathématiques conduit à une indétermination des objets et des énoncés mathématiques. Il s'agit donc d'un raisonnement par l'absurde, l'absurdité de la conclusion (la relativité), devant conduire à rejeter la prémisse (l'adoption des axiomes de Zermelo comme fondement des mathématiques). Skolem peut dès lors conclure :

> « Le principal résultat obtenu ci-dessus est la relativité des notions ensemblistes [*All. : Mengenbegriffe*]. Je l'ai déjà communiqué durant l'hiver 1915-1916 à Göttingen à Mr le Prof. F. Bernstein lors d'une discussion orale. Il y a deux raisons pour lesquelles je ne l'ai pas publié plus tôt : la première est que j'ai été pris entre-temps par d'autres problèmes, la seconde est que je croyais qu'il était suffisamment clair que cette axiomatique ensembliste [*All. : Mengenaxiomatik*] ne pouvait constituer de manière satisfaisante le fondement ultime des mathématiques et que la plupart des mathématiciens ne s'en soucieraient guère. Mais j'ai été surpris de constater ces derniers temps que de très nombreux mathématiciens considéraient ces axiomes de la théorie des ensembles comme le fondement idéal des mathématiques ; il m'a semblé alors qu'il était temps d'en publier une critique. » Skolem 1922, 232.

## L'application à la théorie des ensembles

Un tel usage rhétorique d'un théorème est bien sûr remarquable comme est remarquable cette démonstration de l'*illusion référentielle* en mathématiques, et la certitude de Skolem qu'elle amènera les mathématiciens à renoncer à leurs idées sur la théorie des ensembles. En dépit des enjeux sémiotiques évidents de ces considérations, c'est une autre direction sémiotique que nous voulons suivre pour découvrir à partir de l'application de ce théorème une condition d'un tel énoncé. C'est cette condition que nous allons maintenant mettre en évidence.

Cette condition transparaît déjà au travers d'une critique de Skolem à l'encontre de la manière dont Zermelo a formulé l'un de ses axiomes :

> « Un point très insatisfaisant dans Zermelo est la notion de « proposition définie». Probablement personne ne pourra trouver satisfaisantes les explications données par Zermelo à ce propos. Pour autant que je sache, personne n'a cherché à formuler correctement cette notion, ce qui est très étonnant car cela est très facile à faire et peut, de plus, être fait d'une manière tout à fait naturelle qui se présente d'elle-même. Afin d'expliquer cela - et compte-tenu des considérations précédentes - je mentionne ici les 5 opérations fondamentales de la logique mathématique, pour lesquelles j'utilise les notations de E. Schröder (*Algebra der Logik*). » Skolem 1922

Skolem reproche ici à Zermelo de ne pas avoir défini correctement la notion de « proposition », de s'être contenté d'explications insatisfaisantes. Il oppose qu'il est pourtant facile d'y remédier, tout en notant que cette solution n'a été à sa connaissance reconnue ni par Zermelo, ni par d'autres mathématiciens depuis la publication de ces axiomes. Considérons plus précisément ces axiomes et en particulier celui incriminé pour comprendre à la fois la critique et la solution de Skolem.

Le système de Zermelo comprend sept axiomes : l'axiome d'extensionalité, l'axiome des ensembles élémentaires, l'axiome de séparation, l'axiome de l'ensemble des parties, l'axiome de la réunion, l'axiome du choix et l'axiome de l'infini. La critique porte sur l'axiome de séparation :

> « Si l'énoncé ouvert $\mathcal{E}$ est bien défini pour tous les éléments d'un ensemble $M$, $M$ a toujours un sous-ensemble $M_{\mathcal{E}}$ qui contient tous les éléments $x$ de $M$ pour lesquels $\mathcal{E}(x)$ est vraie, et ceux-là seulement » (Zermelo 1908, trad. 373).

Cet axiome affirme qu'étant donné un ensemble, par exemple celui des nombres entiers, et une proposition, par exemple « être un nombre impair », il existe alors un ensemble constitué des nombres entiers impairs. Il s'agit par ce moyen de ne recevoir que des sous-ensembles d'ensembles déjà formés afin d'éviter la formation d'un ensemble comme l'ensemble de tous les ensembles qui conduit à des antinomies. Avec cet axiome, on ne peut tout au plus que définir l'ensemble de tous les *sous*-ensembles d'un ensemble, c'est-à-dire l'ensemble de ses parties. Skolem reproche donc à Zermelo de ne pas préciser ce que sont ces « énoncés bien définis » auxquels cet axiome se rapporte. Il serait en effet assez peu satisfaisant de fonder les mathématiques sur un système comprenant un axiome qui implique la notion de « proposition définie » sans la caractériser d'aucune manière. Mais il paraît inversement exorbitant de demander une telle caractérisation qui revient ni plus ni moins à une description de tous les énoncés mathématiques passés, présents et à venir. Pourtant, sans cela, ces axiomes ne sauraient prétendre servir de fondement aux mathématiques. Si la caractérisation est trop étroite, se satisfaisant de conditions suffisantes, des ensembles admis par les mathématiciens seraient exclus. Si elle est trop large, se satisfaisant de conditions nécessaires, des ensembles exclus par les mathématiciens seraient admis. Il faut ici une représentation exacte de toutes les propositions mathématiques. Formuler et satisfaire une telle exigence n'apparaît possible qu'en mathématiques. Déterminer les conditions qui y répondent doit nécessairement nous conduire à des conditions qui rendent compte de certaines particularités des énoncés mathématiques. Voici donc la « formulation correcte » qui manquait à Zermelo et que Skolem est en mesure de proposer :

> *« Par une proposition définie nous pouvons maintenant entendre une expression finie construite à partir de propositions élémentaires de la forme $a \varepsilon b$ ou $a = b$ au moyen des 5 opérations mentionnées* [conjonction, disjonction, négation, quantifications universelle et existentielle]*. C'est une notion parfaitement claire et suffisamment générale pour permettre de mener toutes les démonstrations habituelles de théorie des ensembles. J'adopte par conséquent ce point de vue »* Skolem 1922

Skolem affirme donc qu'avec les cinq opérations proposées il est possible d'écrire *tous* les énoncés mathématiques. Avec des notations un peu différentes des siennes, différence sans importance pour notre propos, cela revient à dire que

tous les énoncés mathématiques sont du type :

$$\forall x \forall y (\forall z \, (z \in x \leftrightarrow z \in y) \to x = y) \, ;$$

$$\forall x \forall y \exists z \, (z \in x \wedge y \in z) \, ;$$

$$\exists x \, (0 \in x \wedge \forall y \in x (x \cup \{x\} \in x)) \, ;$$

$$\forall x \exists y \forall z (z \subset x \to z \in y)$$

Et la réciproque doit aussi être vraie pour que les axiomes puissent remplir leur fonction : toutes les propositions que l'on peut écrire de la sorte doivent être des propositions mathématiques bien définies.

Cette critique adressée à la formulation des axiomes par Zermelo permet à Skolem d'introduire dans ceux-ci une autre représentation des propositions mathématiques. Le recours à cette représentation est aussi la condition pour *appliquer* le théorème de Löwenheim-Skolem à ces axiomes, c'est-à-dire la condition que nous voulions découvrir. En effet, il faut pour appliquer le théorème réécrire ces axiomes comme des «équations» ou des « formules » du type de celles qui interviennent dans ses hypothèses, et qui sont aussi celles qui en permettent la démonstration. Cela doit d'ailleurs être fait pour *chacun* des sept axiomes et non seulement pour l'axiome de séparation. Les formules que nous avons données correspondent respectivement aux axiomes d'extensionalité, de la paire (un des ensembles élémentaires), de l'infini et des parties. L'axiome de séparation visé par la critique de Skolem se distingue néanmoins bien des autres puisqu'il implique la *totalité* des énoncés mathématiques. Sa réécriture ne requiert donc pas seulement comme les six autres que l'on sache mettre sous cette forme *un* énoncé mathématique donné mais que l'on sache effectivement le faire pour *tous* les énoncés mathématiques. L'axiome ne va pas être remplacé par un seul axiome, mais par une liste d'axiomes avec *un axiome pour chaque énoncé mathématique*, c'est-à-dire en écrivant la formule suivante pour chaque énoncé mathématique $\varphi$, $\varphi$ devant être effectivement écrit dans la représentation considérée pour que l'axiome puisse lui-même être écrit :

$$\forall A \forall w_1, \ldots, w_n [\forall x \in A \exists ! y \varphi \to \exists Y \forall x \in A \exists y \in Y \varphi]$$

Ainsi, l'application du théorème de Löwenheim-Skolem se fait sous la condition *de disposer d'une représentation de toutes les formules mathématiques*. Cette représentation est la fève dans laquelle tous les énoncés mathématiques peuvent être vus. *Sans elle, le théorème de Löwenheim-Skolem ne pourrait être appliqué aux « axiomes de Zermelo ».* Avec elle les sept axiomes sont transformés en une liste infinie (dénombrable) d'axiomes, à laquelle le théorème de Löwenheim s'applique comme l'a démontré Skolem dans son article précédent.

## Les notions de thèse et d'énoncé inaugural

L'étude de l'application du théorème de Löwenheim-Skolem aux axiomes de la théorie des ensembles a permis de mettre en évidence l'émergence d'une conviction d'un caractère particulier consistant à considérer que l'on dispose d'*une représentation de toutes les propositions mathématiques*. Nous proposons d'appeler une telle conviction une *thèse* et un *énoncé inaugural* l'énoncé qui énonce une thèse. Cette thèse particulière est la *thèse de Frege-Peirce-Schröder-Whitehead-Russell*, ou plus brièvement la *thèse de Russell*. Le terme de *thèse* a été introduit par Kleene pour désigner la *thèse de Church-Turing* qui affirme que *tous les algorithmes* peuvent être représentés par des machines de Turing (resp. fonctions générales récursives, fonctions lambda-définissables, etc.). C'est le même type d'assertion pour les algorithmes au lieu des propositions. Il existe de nombreuses autres thèses comme celles-ci, comme la thèse de Fourier qui propose une représentation de toutes les fonctions par des séries trigonométriques. Dans *La Géométrie* Descartes introduit aussi la thèse selon laquelle toutes les courbes géométriques peuvent être représentées par des équations algébriques. Une autre thèse consiste à considérer que l'on dispose d'une représentation de toutes les démonstrations, thèse que l'on peut appeler *thèse de Frege-Whitehead-Russell*, et que, pour simplifier, nous inclurons ici dans la thèse de Russell (bien qu'il importe par ailleurs de les distinguer). Ces exemples établissent que la logique n'a ici aucun statut particulier et que ce qui a été dégagé à propos de la représentation des propositions mathématiques vaut de la même manière pour les courbes géométriques, les fonctions, les algorithmes, etc. La similarité de ces diverses thèses n'est pas un fait logique, mais un fait sémiotique[4]. Elle n'atteste pas du caractère logique du développement des mathématiques, mais au contraire d'enjeux proprement sémiotiques.

Chacune des représentations associées à ces thèses est comme une fève qui permet aux mathématiciens de considérer de manière *simultanée* et *uniforme* la *totalité* des propositions mathématiques, des courbes géométriques, des algorithmes, etc. On peut entrevoir que quand on dispose de telles représentations il devient possible de produire des énoncés remarquables, démontrables et applicables. La comparaison avec la fève invite à considérer ces représentations comme des images et à les analyser comme telles. Nous proposons de dégager quelques-unes des caractéristiques des thèses qui nous semblent intéressantes en elles-mêmes et susceptibles de servir, au moins à titre de comparaison, pour une analyse de la visualisation en mathématiques. Nous le ferons en nous en tenant principalement à la thèse de Russell pour ne pas multiplier les références mathématiques et historiques, mais la plupart des remarques que nous ferons sont valables pour les autres thèses. Quelques remarques propres à la logique pourront aussi être faites qui peuvent avoir un intérêt sémiotique.

Une thèse est introduite par un *énoncé inaugural*. L'énoncé inaugural de la thèse de Russell a bien été énoncé par Russell :

> "La mathématique pure est la classe de toutes les propositions de la forme "*p* implique *q*", où *p* et *q* sont des propositions contenant une ou plusieurs variables, les mêmes dans les deux propositions, et où ni *p* ni *q* ne contiennent d'autres constantes que des constantes logiques. Et les constantes logiques sont toutes ces notions qui peuvent être définies au moyen de l'implication, de la relation d'un terme à une classe dont il est membre, de la notion de *tel que,* de la notion de relation, et de toutes les autres notions que peut impliquer celle, générale, de proposition de cette forme."
> *Principa Mathematica*, trad. Russell, Bertrand & Roy, Jean-Michel (avant-propos, traduction), *Ecrits de logique philosophique*. Paris, PUF, 1989, p. 21

Ou encore :

---
[4] Nous reprenons ici quelques résultats d'une étude à paraître consacrée aux textes et aux énoncés inauguraux.

> "Au moyen de dix principes de déduction et de dix autres prémisses de nature générale (par exemple, "l'implication est une relation") la totalité de la mathématique peut être rigoureusement et formellement déduite. Et toutes les idées qui figurent en elle peuvent être définies au moyen de celles qui figurent dans ces vingt prémisses. Et ici par mathématique il faut entendre non seulement l'arithmétique et l'analyse, mais aussi la géométrie, euclidienne et non euclidienne, la dynamique rationnelle, et un nombre indéfini d'autres disciplines qui n'ont pas encore vu le jour qui sont dans leur enfance. Le fait que la mathématique n'est dans sa totalité rien qu'autre que la logique symbolique, est une des grandes découvertes de notre temps. Et une fois que cela a été établi, les reste des principes de la mathématique se réduit à l'analyse de la logique symbolique elle-même." *Principia Mathematica*, trad. fr p. 23

Bien que récurrent, ce type d'énoncé échappe aux typologies logiques ou épistémologiques classiques : axiomes, définitions, théorèmes, hypothèses, lois, etc. Il participe de chacune de ces catégories mais ne correspond à aucune exactement. C'est néanmoins un type d'énoncé caractéristique qui accompagne l'introduction de certaines représentations. Cela ressort des termes employés par Skolem. Il reproche à Zermelo des « explications » peu satisfaisantes, demande une « formulation correcte de cette notion » et introduit la représentation qu'il en propose (reprenant celle introduite par Schröder) sous la forme d'un énoncé qui peut être pris pour une définition mais qui a besoin d'être et qui est bien une thèse. Ainsi sa critique peut sembler dénoncer l'absence d'une définition, mais c'est d'un énoncé inaugural dont il a besoin et c'est bien un tel énoncé qu'il introduit sous l'apparence d'une définition. De même, l'énoncé inaugural pourrait être réfuté, mais il ne peut évidemment pas être démontré. Au contraire, il est comme on l'a vu d'une certaine manière une porte d'entrée qui ouvre sur la possibilité de formuler et de démontrer des énoncés mathématiques. Il marque un moment sémiotique fondateur où une représentation mathématique est introduite. Sa caractéristique est d'affirmer la *conformité* de la représentation proposée avec toutes les instances (considérées comme) intuitives d'une notion. Son énoncé affirme que les entités intuitives considérées sont exactement, terme à terme, représentées par la représentation proposée : toute expression de la représentation a un contenu intuitif et toutes les caractéristiques du contenu intuitif peuvent être exprimées de la sorte. C'est donc un énoncé sur la nature de la représentation introduite et ses caractéristiques sémiotiques, avec en particulier la position explicite de *référents*, nullement tenus pour illusoires, de *signifiants*, avec affirmation de leur *conformité*. Les textes qui introduisent de tels énoncés développent un *réalisme* qui pourrait sembler étonnant en mathématiques mais qui est bien attesté et récurrent. Russell peut ainsi déclarer : « *Logic is concerned with the real world just as truly as zoology* »[5]. On peut citer Joseph Fourier pour montrer que ce réalisme n'est propre ni à Russell, ni à la thèse de Russell et qu'il se retrouve aussi bien par exemple à propos des séries trigonométriques près d'un siècle auparavant :

> "L'analyse mathématique a donc des rapports nécessaires avec les phénomènes sensibles ; son objet n'est point créé par l'intelligence de l'homme, il est un élément préexistant de l'ordre universel, et n'a rien de contingent et de fortuit ; il est empreint dans toute la nature. (...) on ne pourrait apporter aucun changement dans la forme de nos solutions, sans leur faire perdre leur caractère essentiel, qui est de représenter les phénomènes." Fourier, Joseph, *Théorie de la chaleur* (2ème

---
5   Russell, Bertrand, *Introduction to mathematical philosophy*. London, New York : Allen and Unwin, Macmilan, 1919, 169.

Il importe, tant pour l'étude des thèses que pour l'étude des conditions de possibilité des énoncés généraux en mathématiques que nous considérons ici de distinguer le moment fondateur, celui où la thèse est énoncée pour la première fois, de la reprise qui pourra ensuite être faite de la représentation que l'énoncé inaugural a servi à introduire et de ses conséquences. La thèse a en effet d'abord besoin d'être soutenue. Les trois épais volumes des *Principia Mathematica* de Whitehead & Russell sont ainsi largement consacrés à soutenir leur thèse, c'est-à-dire à reproduire, pas à pas, une à une, toutes les propositions mathématiques et leur démonstration. Cela donne lieu à un type de texte particulier, les *textes inauguraux*. Indépendamment des démonstrations reproduites, les *Principia* sont surtout eux-mêmes la preuve de ce qu'ils défendent. Ils doivent être *figuratifs* pour montrer qu'ils peuvent reproduire exactement, sans manques ni ajouts, ce qu'ils représentent. Ainsi, attestent-ils qu'il est *effectivement* possible d'écrire toutes les mathématiques, comme nous l'avons fait nous-mêmes pour quelques-uns des axiomes de Zermelo.

*Exemple de texte figuratif. Pages extraites des* Principia Mathematica *de Whitehead & Russell.*

C'est une preuve éminemment singulière qui n'a pas besoin et qui ne saurait être vraiment répétée, et qui n'est pas non plus ouverte à la généralisation. Mais surtout, une fois donnée, elle perd de son intérêt, et avec elle le texte qui la présente. L'accomplissement de sa tache lui fait perdre sa fonction qui est autant son sens que sa raison d'être. Ces textes ont ainsi un caractère performatif dont l'action réside dans la charge qu'ils confèrent à la représentation qu'ils introduisent en établissant sa conformité et en montrant quelques-uns des bénéfices que l'on tire de son usage, propre à chaque thèse et à chaque représentation. Leur performance ne tient pas à des conventions préalables, ce qui ne les exclut pas, et consiste à faire le nécessaire pour l'introduction d'une nouvelle représentation mathématique. Ils ouvrent une porte qui, n'ayant pas ensuite à être fermée, n'apparaîtra pas comme ouverte, ni comme ayant été ouverte et ayant dû l'être. Ils accomplissent malgré eux l'idéal du *hacker*, mais qui n'est que rarement celui attribué au scientifique : pénétrer et transformer un système sans y laisser de trace. C'est le destin tragique des textes inauguraux qui instaurent des représentations qui seront ensuite adoptées, au travers desquelles ils seront ensuite eux-mêmes relus pour de ce fait sembler ne plus rien faire ou faire en trop de mots ce qu'il est possible de faire en peu mais en fait en tirant parti de ce que la représentation qu'ils ont introduit synthétise.

Les textes de Löwenheim et de Skolem, ou encore ceux de Post ou de Gödel où sont démontrés divers théorèmes de complétude ou d'incomplétude sont d'un tout autre genre. Ils peuvent d'emblée déclarer : *voilà toutes les propositions!* Ils bénéficient des possibilités offertes par la représentation devenue disponible. Skolem peut, comme nous l'avons vu, considérer en 1922 qu'il est « très facile » de représenter toutes les propositions mathématiques, que cela peut « être fait d'une manière tout à fait naturelle qui se présente d'elle-même », et reprocher à Zermelo de ne pas l'avoir fait en 1908. De la même manière Löwenheim pouvait en 1915 simplement affirmer dans l'article où il énonce son théorème :

> "Il semble que toutes les questions importantes des mathématiques et du calcul logique peuvent être ramenées à de telles équations relationnelles." Löwenheim 1915

Emil Post pourra de la même manière donner à voir toutes les propositions[6] :

> « Il est souhaitable dans ce qui suit d'avoir devant nous la vision de la totalité de ces fonctions jaillissant d'une variable inchangée $p$ en des formes d'une complexité toujours croissante pour constituer le tableau triangulaire infini
>
> $$p$$
>
> $$p \vee p,\ p_1 \vee p_2,\ \sim p$$
>
> $$p \vee \sim p,\ ...,\ \sim p_1 \vee \sim p_2,...,\ (p_1 \vee p_2) \vee (p_3 \vee p_4),\ \sim(p_1 \vee p_2),\ \sim(p \vee p),\ \sim\sim p$$
>
> ..................................................... »
>
> Post 1921

L'article dans lequel Gödel démontre en 1931 ses célèbres théorèmes d'incomplétude peut aussi commencer ainsi :

> « Le développement des mathématiques vers plus d'exactitude a conduit, comme nous le savons, à en formaliser de larges secteurs, de telle sorte que la démonstration puisse s'y effectuer uniquement au moyen de quelques règles mécaniques. Les systèmes formels les plus complets établis jusqu'à ce jour sont, d'un côté, le système des *Principia Mathematica* (*PM*) et, de l'autre, le système axiomatique de la théorie des ensembles établi par Zermelo-Fraenkel (et développé par J. von Neumann). » Gödel 1931, trad. 107

---

6  Il se restreint en fait ici à une partie des propositions (appelée calcul propositionnel), mais cette *restriction* atteste aussi bien de la réception de la représentation utilisée, qui et en l'occurrence celle des *Principia,* et Post entendait à ce moment là étendre son théorème à toutes les propositions des *Principia.*

Ces citations attestent autant de la réception des représentations associées à la thèse de Russell que de l'oubli des conditions de son introduction. Et les théorèmes énoncés et démontrés dans chacun d'eux exploitent les possibilités offertes par le fait de considérer avoir une représentation de toutes les propositions et, pour les théorèmes de Gödel, de toutes les démonstrations. Et l'on ne saurait trop souligner que ces mathématiciens peuvent avoir à ces représentations, et *à ce qu'elles représentent*, un rapport que ceux qui les ont introduites auraient difficilement pu d'emblée avoir. Whitehead & Russell ne pouvaient disposer de la représentation qu'ils soutenaient. L'introduction d'une telle représentation est un processus. Et si les *Principia Mathematica* ou l'*Algebra der Logik* perdent une partie de leur sens, inversement, les théorèmes qui exploitent leurs représentations n'y ont pas non plus leur place. Ces théorèmes ne sont pas comparables à ceux qui y sont démontrés, pourtant nombreux et censés les recouvrir tous. Ce ne sont pas des maillons supplémentaires de la chaîne déductive déroulée dans les *Principia* : ce sont pour une part des conséquences de la représentation qui y est soutenue. Reconnaître cette différence c'est en particulier reconnaître que les mathématiques ne sont pas symboliques. Conséquence que renforcent les différences entre les représentations introduites par les divers énoncés inauguraux qu'une conception symbolique est amenée à ignorer. Le choix d'un théorème de logique présente ici un intérêt particulier. En effet, les *Principia Mathematica* constituent sans doute l'un des textes le plus conforme à une conception symbolique des mathématiques. Mais les *théorèmes* qui illustrent par exemple la généralité particulière des énoncés mathématiques, comme celui de Löwenheim-Skolem, ne sont pas de ceux que l'on trouve dans ce texte et les moyens employés pour les démontrer ne sont pas non plus les mêmes. Ils disposent de représentations dont ne disposaient pas de la même manière les auteurs des *Principia*. La critique que Skolem adresse à Zermelo permet aussi de mieux comprendre en quoi la représentation associée à la thèse de Russell intervient comme condition de possibilité de théorèmes et de préciser l'enjeu des changements sémiotiques dans et pour l'histoire des mathématiques. L'article de Zermelo est sans conteste un texte mathématique. Il énonce et démontre divers théorèmes dont on pourrait d'ailleurs aussi chercher à déterminer leurs conditions de possibilité sémiotiques. La reformulation de ses axiomes n'introduit donc évidemment pas les conditions sémiotiques à partir desquelles la formulation et la démonstration d'énoncés généraux deviendrait enfin possible ; elle ne fait pas basculer le texte de l'extérieur vers l'intérieur des mathématiques. En revanche, elle modifie bien comme on l'a vu les possibilités d'énonciation et de démonstration. Le fait qu'il devient dès lors possible, du fait de cette réécriture, de leur appliquer le théorème de Löwenheim-Skolem suffit à l'établir. Bien d'autres développements mathématiques pourraient être donnés en exemples. On voit ainsi de quelle manière des changements sémiotiques interviennent en mathématiques, contribuent à leur évolution et donc participent de leur histoire, en même temps qu'ils passent inaperçus et sont à peu près systématiquement ignorés. En l'occurrence, le fait de (considérer) disposer d'une représentation de toutes les propositions mathématiques permet à Skolem de modifier radicalement le statut des axiomes de Zermelo par une réécriture que sa critique adressée à Zermelo permet de rendre visible. Cette opération sera ensuite encore plus transparente et ses conditions et ses enjeux systématiquement ignorés. La lecture de l'article de Zermelo, et l'évidence aveuglante de la formulation de ses axiomes, ne suffit pas pour découvrir la nécessité d'une réécriture que le lecteur actuel aura tendance à opérer spontanément. Il est pour cela utile de considérer un texte comme celui de Skolem pour saisir, dans le reflet d'une critique, une représentation déjà devenue transparente pour son auteur.

L'exemple du théorème de Löwenheim-Skolem a permis de mettre en évidence l'enjeu et le rôle de la représentation d'une *totalité*, en l'occurrence des propositions mathématiques, dans la possibilité de produire et de démontrer des énoncés d'une généralité qui paraît assez propre aux mathématiques et de dégager quelques éléments de son historicité. Ces représentations n'épuisent pas l'analyse, même sémiotique, des formes d'expression de la généralité en mathématiques. Elles permettent néanmoins de repérer une caractéristique qui se retrouve dans des parties suffisamment variées des mathématiques associée à de nombreux moments fondateurs de leur histoire. Ces représentations ne sont pas le produit d'une analyse de sémioticien. Elles témoignent au contraire d'un enjeu régulièrement thématisé dans les textes mathématiques et repérable par un énoncé inaugural d'un type remarquable justiciable d'une analyse pragmatique. Mais leur évolution indique aussi un oubli tout aussi systématique de leur introduction et de ses conditions, oubli qui se manifeste aussi bien dans les développements mathématiques auxquels elles contribuent que dans les travaux d'histoire des mathématiques qui en rendent compte. L'intervention du sémioticien a ici surtout consisté à ne pas se laisser abuser par la transparence de ces représentations et à en saisir quelques reflets. Loin de confirmer l'idée d'une mathématique symbolique, les thèses et les énoncés inauguraux montrent au contraire l'insuffisance d'une telle conception. Les représentations introduites ne correspondent pas en effet au modèle d'un code, comme un alphabet en serait idéalement un exemple en raison notamment de leur signification d'ensemble à laquelle ne prétend pas un code. Elles ne prétendent pas seulement par exemple permettre la transcription de tous les sons d'une langue, voire de toutes les langues ; les énoncés inauguraux affirment aussi la condition réciproque (et qui est d'ailleurs la partie qui pose au premier abord le moins de problèmes) : tous les sons que l'on peut composer avec la représentation introduite sont effectivement des sons de la langue ou d'une langue. Autrement dit, les énoncés inauguraux introduisent des représentations qui n'excèdent pas ce qu'elles servent à représenter. Ce sont des représentations qui peuvent pour cette raison être tenues pour sans effets, et sembler pour cela totalement et intégralement réalistes. Ce n'est qu'en étant reprises, et de ce fait transformées, qu'elles perdront une part de cette conformité, ce qui ne doit pas pour autant faire oublier que celle-ci aura été la raison de leur introduction. La récurrence tout au long de l'histoire de ces moments réalistes fondateurs contredit aussi l'idée d'un développement des mathématiques suivant une abstraction croissante, qui n'est jamais que le moyen d'ignorer l'histoire des mathématiques et ses conditions sémiotiques.

## Enoncé inaugural et théorème de représentation

Revenons maintenant à l'énoncé euclidien qui affirme que « toute figure rectiligne peut être transformée en carré ». Si sa démonstration ne se fait pas par induction, ou par un raisonnement par l'absurde, on peut chercher la *représentation* qui sert de fève et qui permet de considérer *toutes* les figures. La démonstration de cet énoncé ne saurait créer de la généralité et, sauf à croire à quelques effets mystérieux, sur lesquels toute une philosophie des mathématiques pourra ensuite se déployer, la généralité doit y être introduite par une représentation. Quand on sait qu'une telle fève peut exister, il n'est plus difficile de la trouver. La démonstration, dans les *Eléments* d'Euclide, de ce théorème repose sur la possibilité de *décomposer toute figure rectiligne en*

*triangles*. Mais l'énoncé disant que « toute figure rectiligne peut être décomposée en triangles » n'est pas un énoncé inaugural. La représentation d'une figure rectiligne par des triangles n'est pas une nouvelle représentation des figures rectilignes comme les formules des *Principia* le sont des propositions mathématiques. Les formules logiques sont exprimées dans un *autre* système d'expressions que ne le sont les propositions mathématiques (ce qui était plus évident encore avant qu'on ne se mette à les confondre...). Les figures rectilignes décomposées en triangles sont au contraire constituées des *mêmes expressions* que les figures rectilignes. Il ne s'agit pas cette fois d'un énoncé inaugural mais d'un *théorème de représentation*. La fève ne se trouve donc pas là. Elle se trouve en l'occurrence dans la représentation des figures rectilignes au moyen de segments de droites juxtaposés par leurs extrémités. Le mathématicien dispose en effet ainsi d'une représentation de *toutes* les figures rectilignes qui rend possible la formulation et la démonstration d'énoncés généraux, c'est-à-dire de théorèmes.

## II. L'expression comme problème

Les thèses et les énoncés inauguraux ainsi que les théorèmes de représentation ont permis de mettre en évidence le rôle joué par des représentations de *totalités* dans l'expression de la généralité en mathématiques. Le rôle de ces représentations se comprend peut-être mieux si on les considère comme des *expressions* d'une totalité, comme la fève en est une d'un paysage. Ainsi, les *Principia Mathematica* introduisent des expressions pour les propositions et les démonstrations mathématiques qui vont ensuite pouvoir être considérées pour elles-mêmes et servir d'expression de la totalité des propositions et des démonstrations mathématiques. Leur intérêt réside dans cette représentation totale et dans ses caractéristiques propres qu'il appartient à la sémiotique de définir et de dégager. Une part de la créativité du mathématicien est de les avoir introduites. De la même manière, les équations algébriques de Descartes sont une expression des courbes géométriques, les séries trigonométriques des fonctions, etc. On peut ainsi traiter de manière symétrique les deux termes de la thèse en les considérant l'un et l'autre comme un système d'expressions avec leurs caractéristiques et leur intérêt propre. Les énoncés inauguraux sont dès lors des énoncés qui tout au long de l'histoire des mathématiques accompagnent l'introduction d'un nouveau *système d'expressions* et qui offre la possibilité de nouveaux développements. Réinscrits dans une problématique générale de l'*expression,* ils se distinguent par le fait d'être associés à l'introduction de *systèmes* d'expressions, c'est-à-dire, d'une pluralité d'expressions constituant d'une certaine manière elles-mêmes une totalité autonome, qui peut être considérée pour elle-même, sans doute pourvue d'une syntagmatique et d'une paradigmatique. Ces systèmes déterminent à la fois le pouvoir d'expression et le pouvoir de résolution des mathématiques. Mais les enjeux de l'expression ne se réduisent pas à celui de l'introduction d'un *système*. C'est ce que nous voudrions indiquer en rappelant le résultat d'une analyse sémiotique des *Disquisitiones arithmeticae* de Gauss[7].

---

7  Voir Alain Herreman, « Vers une analyse sémiotique de la théorie des ensembles : hiérarchies et réflexivité ». *Philosophia Scientia* 9 (2), 165-187, 2005 (http://perso.univ-rennes1.fr/alain.herreman/hte.pdf).

Dans ce livre, fondateur de nombreux développements mathématiques du 19ème siècle[8], la généralité est systématiquement exprimée au moyen d'expressions du type *2n+1*. Ces expressions ont plusieurs caractéristiques remarquables. Elles expriment une totalité, ici la totalité des nombres impairs. Ainsi, une expression est associée à une infinité d'autres : celles des nombres impairs. Ce que l'expression « nombre impair » fait aussi très bien. Mais en substituant au *n* de *2n+1* toutes les expressions pour les nombres on obtient toutes les expressions des nombres impairs (on notera l'intervention de l'énoncé inaugural affirmant que l'on dispose d'une expression pour tous les nombres entiers) ; l'expression *2n+1* est une expression d'invariante au moyen de laquelle on peut dériver par substitution toutes les variétés dont elle est l'invariante. Cette fois, l'expression « nombre impair » ne fait pas cela. Ce n'est pas tout : *2n+1* est aussi un nombre impair au sens où il est aussi possible de faire avec cette expression toutes les opérations arithmétiques qu'il est possible de faire avec un nombre impair particulier sans faire intervenir sa valeur particulière. Autrement dit, l'expression de l'invariante est telle que l'on peut lui attribuer la propriété dont elle est l'expression d'invariante. Ce n'est à nouveau pas le cas de l'expression « nombre impair » qui n'est pas un nombre impair. L'expression d'invariante génératrice est un nombre impair sans être aucun nombre impair particulier. Ainsi, *2n+1* peut à la fois représenter les nombres impairs et en avoir la propriété, mais sans en faire partie. Ces expressions peuvent être comparées aux expressions d'invariante représentantes dont elles sont proches tout en s'en distinguant nettement. Ce dernier type d'expression est aussi souvent utilisé, y compris en mathématiques, et consiste à prendre l'expression d'un nombre impair, par exemple « 9 » ou « neuf », mais à ne faire intervenir que les propriétés communes à tous les nombres impairs, en faisant abstraction de sa valeur particulière (on peut parler ailleurs de *specimen*, de prototype, de modèle etc.). Une variété sert ici d'invariante, et donc l'invariante est aussi dans ce cas une variété comme pour l'expression *2n+1* à ceci près que celle-ci a une expression propre et qu'elle a la propriété d'engendrement des variétés par substitution indiquée plus haut.

Bien sûr, les expressions considérées par Gauss sont bien plus complexes, mais ce sont bien là leurs caractéristiques sémiotiques communes. Néanmoins, il lui arrive aussi d'avoir des problèmes d'expression… Il lui arrive par exemple d'avoir à considérer des familles de nombres pour lesquelles il ne dispose pas d'expression de ce type. Si par exemple au lieu de considérer seulement la propriété « être un nombre impair » on ajoute aussi celle d'« être un carré » pour considérer la propriété « être un nombre impair et être un carré », on a encore dans ce cas une expression d'invariante génératrice : *(2n+1)²*. Il se trouve en effet que cette expression donne exactement tous les nombres carrés impairs. Cela tient en fait aux valeurs particulières « 2 » et « 1 ». Si on remplace maintenant celles-ci par deux nombres quelconques *x* et *y* et que l'on souhaite à nouveau trouver une expression unique pour les nombres de la forme *xn+y* qui soient aussi des carrés, alors l'expression *(xn+y)²* ne convient plus. Que fait Gauss dans ce cas ? Il considère l'*ensemble* des nombres carrés de la forme *xn+y* et le désigne aussi par une lettre, par exemple *ω*. Cette lettre, considérée comme un symbole, est un symbole comme un autre, mais elle n'a plus du tout les remarquables propriétés sémiotiques des expressions d'invariante génératrices. Son introduction traduit une sorte de dépit sémiotique. Ainsi, Gauss utilise un mode d'expression dont il a hérité, mais l'usage qu'il en fait l'amène à avoir un *problème d'expression*, c'est-à-dire que le mode d'expression

---

reçu, qu'il exploite magistralement, ne lui fournit pas toutes les expressions dont il a besoin. Il y a des éléments du paysage qui n'entrent pas dans la fève. Gauss doit renoncer ou trouver des expédients sémiotiques et, *bien malgré lui*, recourir à d'autres expressions. Il recourt alors *systématiquement* à des *ensembles*. Ainsi, les ensembles s'introduisent dans les *Disquisitiones*, c'est-à-dire dès 1801, comme la solution systématique aux problèmes d'expression rencontrés par Gauss dans l'usage de son système d'expression. Ce recourt ne s'accompagne d'aucun énoncé inaugural. Il n'est pas thématisé, il reste subreptice et les ensembles ne constituent pas ici un *système*. En particulier, ils ne définissent ni une syntagmatique, ni une paradigmatique et ne conduisent pas à considérer les problèmes qui le sont quand c'est le cas. Compte-tenu de l'importance souvent accordée à l'infini dans le développement de la théorie des ensembles considéré à partir des travaux de Cantor, il convient peut-être de remarquer que les raisons de recourir à des ensembles sont en l'occurrence indépendantes de leur caractère fini ou infini : le problème d'expression peut se poser et se pose effectivement aussi avec des ensembles finis (il suffit de considérer les nombres entiers carrés de la forme *xn+y et plus petits*, par exemple, *que $(x+y)^2$*). On observera enfin à nouveau qu'une conception symbolique des mathématiques conduit encore à ignorer à la fois les différents types d'expressions que nous avons distingués et le rôle dans l'histoire des mathématiques des problèmes d'expression auxquels ils donnent lieu ou qu'ils résolvent, avant d'en poser de nouveaux..., et qui sont autant d'obstacles et de tremplins sémiotiques objectifs et repérables.

D'autres expressions de la généralité en mathématiques ont été diversement décrites et restent à décrire. Mais il appartient aussi à la sémiotique de les décrire qu'elles aient ou non été introduites par des énoncés inauguraux. Nous nous sommes délibérément contentés ici du terme générique de « représentation », et nous n'avons en quelque sorte fait que montrer celles-ci (ou leur description) du doigt. Elles doivent bien sûr être décrites, caractérisées pour ainsi mieux mettre en évidence leurs différences et le rôle spécifique de chacune dans le développement des mathématiques.